\documentclass[12pt]{amsart}
\usepackage{amsmath}
\usepackage{amssymb}
\usepackage{stmaryrd}
\usepackage{epsfig,color,esint}
\usepackage{blindtext}
\usepackage{enumerate}
\usepackage{enumitem}
\usepackage{hyperref}
\usepackage{graphicx}
\usepackage{url}
\usepackage{bbm}
\usepackage{nicefrac,mathtools}
\usepackage{mathrsfs}

\usepackage{comment}
\usepackage{pgfplots}
\usepackage{tikz}
\usetikzlibrary{arrows,shapes,trees,backgrounds}

\numberwithin{equation}{section}

\newtheorem{theo}{Theorem}[section]

\newtheorem{coro}{Corollary}[section]

\def\begeq{\begin{equation}}
\def\endeq{\end{equation}}

\begin{document}

\title{Non-quadratic Euclidean complete affine maximal type hypersurfaces for $\theta\in(0,(N-1)/N]$}
\author{Shi-Zhong Du}
\thanks{The author is partially supported by NSFC (12171299), and GDNSF (2019A1515010605)}
  \address{The Department of Mathematics,
            Shantou University, Shantou, 515063, P. R. China.} \email{szdu@stu.edu.cn}

\renewcommand{\subjclassname}{%
  \textup{2010} Mathematics Subject Classification}
\subjclass[2010]{53A15 $\cdot$ 53A10 $\cdot$ 35J60}
\date{Jul. 2022}
\keywords{Affine maximal hypersurfaces, Bernstein problem.}

\begin{abstract}
   Bernstein problem for affine maximal type equation
     \begin{equation}\label{e0.1}
      u^{ij}D_{ij}w=0, \ \ w\equiv[\det D^2u]^{-\theta},\ \ \forall x\in\Omega\subset{\mathbb{R}}^N
     \end{equation}
  has been a core problem in affine geometry. A conjecture proposed firstly by Chern (Proc. Japan-United States Sem., Tokyo, 1977, 17-30) for entire graph and then extended by Trudinger-Wang (Invent. Math., {\bf140}, 2000, 399-422) to its full generality asserts that any Euclidean complete, affine maximal type, locally uniformly convex $C^4$-hypersurface in ${\mathbb{R}}^{N+1}$ must be an elliptic paraboloid. At the same time, this conjecture was solved completely by Trudinger-Wang for dimension $N=2$ and $\theta=3/4$, and later extended by Jia-Li (Results Math., {\bf56} 2009, 109-139) to $N=2, \theta\in(3/4,1]$ (see also Zhou (Calc. Var. PDEs., {\bf43} 2012, 25-44) for a different proof). On the past twenty years, much efforts were done toward higher dimensional issues but not really successful yet, even for the case of dimension $N=3$. Recently, counter examples were found in \cite{Du2} (J. Differential Equations, {\bf269} (2020), 7429-7469) for $N\geq3$ and $\theta\in(1/2,(N-1)/N)$ using a much more complicated argument. In this paper, we will construct explicitly various new Euclidean complete affine maximal type hypersurfaces which are not elliptic paraboloid for the improved range
     $$
      N\geq2, \ \ \theta\in(0,(N-1)/N].
     $$
  .
\end{abstract}

\maketitle\markboth{Affine Maximal Hypersurfaces}{Bernstein Theorem}

\tableofcontents

\section{Introduction}

In this paper, we study the local uniformly convex solution to affine maximal type equation
   \begin{equation}\label{e1.1}
       D_{ij}(U^{ij}w)=0, \ \ \forall x\in\Omega\subset{\mathbb{R}}^N,
   \end{equation}
where $U^{ij}$ is the co-factor matrix of $u_{ij}$ and $w\equiv[\det D^2u]^{-\theta}, \theta>0$. Equation \eqref{e1.1} is the Euler-Lagrange equation of the affine area functional
   \begin{eqnarray*}
     {\mathcal{A}}(u,\Omega)&\equiv&\int_\Omega[\det D^2u]^{1-\theta}\\
      &=&\int_{{\mathcal{M}}_\Omega}K_0^{1-\theta}(1+|Du|^2)^{\vartheta}dV_{g_0},\ \ \vartheta=\frac{N+1}{2}-\frac{N+2}{2}\theta
   \end{eqnarray*}
for $\theta\not=1$ and
    $$
     {\mathcal{A}}(u,\Omega)\equiv\int_\Omega\log\det D^2u
    $$
for $\theta=1$, where $g_0$ and $K_0$ are the induced metric and the Gauss curvature of the graph ${\mathcal{M}}_\Omega\equiv\Big\{(x,z)\in{\mathbb{R}}^{N+1}|\ z=u(x), x\in\Omega\Big\}$ respectively. Noting that
   $$
     D_jU^{ij}=0, \ \ \forall i=1,2,\cdots, N,
   $$
equation \eqref{e1.1} can also be written by
   $$
      U^{ij}D_{ij}w=0,
   $$
or equivalent to
   \begin{equation}\label{e1.2}
      u^{ij}D_{ij}w=0
   \end{equation}
for $[u^{ij}]$ denoting the inverse of usual metric $[u_{ij}]$ of graph ${\mathcal{M}}_\Omega$.

 The classical affine maximal case $\theta\equiv\frac{N+1}{N+2}$ has been studied extensively on the past. If one introduces the affine metric
   $$
     A_{ij}=\frac{u_{ij}}{[\det D^2u]^{1/(N+2)}}
   $$
on  ${\mathcal{M}}_\Omega$ comparing to the Calabi's metric $g_{ij}=u_{ij}$ and sets
  $$
    H\equiv[\det D^2u]^{-1/(N+2)},
  $$
it is not hard to see that \eqref{e1.2} turns to be
   \begin{equation}\label{e1.3}
    \triangle_{{\mathcal{M}}}H=0
   \end{equation}
for Laplace-Beltrami operator
   $$
     \triangle_{{\mathcal{M}}}\equiv\frac{1}{\sqrt{A}}D_i(\sqrt{A}A^{ij}D_j)=HD_i(H^{-2}u^{ij}D_j)
   $$
with respect to this affine metric, where $A$ is the determinant of $[A_{ij}]$ and $[A^{ij}]$ stands for the inverse of $[A_{ij}]$. So, the hypersurface ${\mathcal{M}}$ is affine maximal if and only if $H$ is harmonic on ${\mathcal{M}}_\Omega$.

A conjecture proposed by Chern \cite{Ch} for dimension $N=2$ and $\theta=3/4$ asserts that every locally convex entire graph of \eqref{e1.2} must be a paraboloid. Much efforts were done toward this conjecture, (see examples \cite{B,Ca1,Ca2,JL,TW} for partial results) until a landmark paper by Trudinger-Wang \cite{TW}. At there, they strengthened the Chern's conjecture to its full generality  and then give a proof for 2 dimensional affine maximal case. Before our discussion, let us first reformulate this version of full Bernstein problem by Trudinger-Wang to all dimensions $N$ and positive $\theta$ as following.\\

\noindent\textbf{Full Bernstein Problem in Sense of Trudinger-Wang:} Given dimension $N\geq1$ and $\theta>0$, whether any locally uniform convex, Euclidean complete affine maximal type hypersurfaces must be an elliptic paraboloid?\\

Throughout this paper, a hypersurface on ${\mathbb{R}}^{N+1}$ is called to be Euclidean complete if and only if it is complete under the induced metric from standard metric of the whole Euclidean space. At beginning, we first sum some positive results to this Bernstein problem under below, which mainly owed to Trudinger-Wang \cite{TW} for $N=2, \theta=3/4$. (see also \cite{JL2,Z} for $N=2, \theta\in(3/4,1]$)\\

  \noindent\textbf{Theorem A.} Full Bernstein theorem in sense of Trudinger-Wang holds under either one of the following conditions:

  (1) $N=1$ and $\theta>0$, or

  (2) $N=2$ and $\theta\in[3/4,1]$.\\

Unlike the two dimensional case, attempts toward full Bernstein theorem on dimension $N\geq3$ were not yet successful. (see for examples \cite{JL,JL1,JL2,L,Wang,Z}) A first well known partial counter example on dimension $N\geq10$ and $\theta=11/12$ was found by Trudinger-Wang in \cite{TW}, which is defined by the graph of
   $$
     u(y,t)=\sqrt{|y'|^9+t^2}+|\widetilde{y}|^2, \ \ y\equiv(y',\widetilde{y}), y'\in{\mathbb{R}}^9, \widetilde{y}\in{\mathbb{R}}^{N-10}, t\in{\mathbb{R}}.
   $$
Because that this solution is not smooth at origin and has degenerate convexity along the set
   $$
    \Big\{(y,t)\in{\mathbb{R}}^N\Big|\ y=(y',\widetilde{y}), y'=0\Big\},
   $$
we looked for regular counter examples to the full Bernstein problem and obtained the following result in \cite{Du2}.\\

\noindent\textbf{Theorem B.} For each dimensions $N\geq3$, supposing that
    \begin{equation}\label{e1.4}
      \theta\in(1/2, (N-1)/N),
    \end{equation}
there exist Euclidean complete affine maximal type hypersurfaces which are not elliptic paraboloid.\\

 The proof of Theorem B was built on a much more complicated technic and did not produce explicit expressions for the solutions. In this paper, we shall improve the above range of $\theta$ and construct the desired hypersurfaces explicitly to show the following result.\\

\noindent\textbf{Theorem C.} Assuming that $N\geq2$ and
    \begin{equation}\label{e1.5}
      \theta\in(0, (N-1)/N],
    \end{equation}
there exist Euclidean complete affine maximal type hypersurfaces which are not elliptic paraboloid.\\

The proof of Theorem C will be decomposed into several parts. At first, by looking for solutions of Warren type \cite{W}, we handle the case $\theta\in(0,1/2)$ in Section 2 and Section 3. Utilizing the functions of Trudinger-Wang type \cite{TW}, we construct Euclidean complete hypersurfaces for $N\geq2, \theta=1/2 \mbox{ or } (N-1)/N$ in Section 4. Finally, we complete the proof of Theorem C for $\theta\in(1/2, (N-1)/N)$ in the last section. On the range $\theta\in(1/2, (N-1)/N)$, we construct explicitly new non-quadratic affine maximal type hypersurfaces with a completely different argument. As a corollary of Theorem C, one can estimate the critical dimension $N^*(\theta)$ defined in \cite{Du2} as following.\\

\noindent\textbf{Corollary D.} For each $\theta\in(0,1)$, the critical dimension $N^*(\theta)$ satisfies that
  \begin{equation}\label{e1.6}
     \Bigg[\frac{1}{1-\theta}\Bigg]_*-1\geq N^*(\theta)\geq\begin{cases}
       1, & \theta\in(0,3/4)\\
       2, & \theta\in[3/4,1),
     \end{cases}
  \end{equation}
where $[z]_*$ stands for the smallest integer no less than $z$.\\

This corollary follows directly from Theorem A and Theorem C.

\vspace{40pt}

\section{Warren type solutions}

From here now, we turn to construct non-quadratic affine maximal type hypersurfaces. Considering the functions of the form
    $$
     u(y,t)=|y|^2\eta(t)+\varphi(t), \ \ \forall y\in{\mathbb{R}}^{n}, t\in{\mathbb{R}}
    $$
and set
    $$
      x\equiv(y,t)\in{\mathbb{R}}^N, \ \ \ r\equiv|y|, \ \ \ N\equiv n+1
    $$
for short. Solutions in a similar form have been used by Warren in \cite{W} for $k-$Hessian equations, where counter examples of Bernstein property were constructed without convexity by choosing
   $$
     \eta(r)=e^r, \ \ \varphi(t)=\iint\frac{1-Be^{kt}}{Ae^{(k-1)t}}
   $$
for some determined constants  $A, B$.  For affine maximal type equation \eqref{e1.1}, we extend these functions to a more general form in order constructing non-quadratic Euclidean complete hypersurfaces. Direct calculation shows that
  \begin{eqnarray*}
    D^2u&=&\left(
     \begin{array}{ccccc}
       2\eta & 0 & \cdots & 0 & 2y_1\eta' \\
       0  &  2\eta & \cdots & 0 & 2y_2\eta'\\
       \cdots & \cdots& \cdots & \cdots & \cdots\\
       0 & 0 & \cdots & 2\eta & 2y_n\eta'\\
       2y_1\eta' & 2y_2\eta' & \cdots & 2y_n\eta' & |y|^2\eta''+\varphi''
     \end{array}
    \right)\\
    &=&2\eta\left(
     \begin{array}{ccccc}
       1 & 0 & \cdots & 0 & y_1\eta'/\eta \\
       0  &  1 & \cdots & 0 & y_2\eta'/\eta\\
       \cdots & \cdots& \cdots & \cdots & \cdots\\
       0 & 0 & \cdots & 1 & y_n\eta'/\eta\\
       y_1\eta'/\eta & y_2\eta'/\eta & \cdots & y_n\eta'/\eta & \frac{|y|^2\eta''+\varphi''}{2\eta}
     \end{array}
    \right).
  \end{eqnarray*}
So, there holds
   \begin{equation}\label{e7.1}
     \det(D^2u)=(2\eta)^n\Bigg[r^2\Bigg(\eta''-\frac{2\eta'^2}{\eta}\Bigg)+\varphi''\Bigg].
   \end{equation}
By rotating the coordinates if necessary, one may assume that $y=(y_1,0,\cdots,0)$, and then simplify the Hessian matrix $D^2u$ by
  \begin{equation}\label{e7.2}
    D^2u=2\eta\left(
     \begin{array}{ccccc}
       1 & 0 & \cdots & 0 & r\eta'/\eta \\
       0  &  1 & \cdots & 0 & 0\\
       \cdots & \cdots& \cdots & \cdots & \cdots\\
       0 & 0 & \cdots & 1 & 0\\
       r\eta'/\eta & 0 & \cdots & 0 & \frac{r^2\eta''+\varphi''}{2\eta}
     \end{array}
    \right).
  \end{equation}
Thus, the solution is strictly convex if and only if
  \begin{equation}\label{e7.3}
    \eta>0, \ \ \ \eta''-\frac{2\eta'^2}{\eta}\geq0, \ \ \ \varphi''>0.
  \end{equation}
Moreover, the inverse matrix of $D^2u$ is given by
   \begin{equation}\label{e7.4}
    [u^{ij}]=(2\eta)^{-1}\left(
     \begin{array}{ccccc}
       \frac{r^2\eta''+\varphi''}{r^2\eta''+\varphi''-2r^2\eta'^2/\eta} & 0 & \cdots & 0 & -\frac{2r\eta'}{r^2\eta''+\varphi''-2r^2\eta'^2/\eta} \\
       0  &  1 & \cdots & 0 & 0\\
       \cdots & \cdots& \cdots & \cdots & \cdots\\
       0 & 0 & \cdots & 1 & 0\\
       -\frac{2r\eta'}{r^2\eta''+\varphi''-2r^2\eta'^2/\eta} & 0 & \cdots & 0 & \frac{2\eta}{r^2\eta''+\varphi''-2r^2\eta'^2/\eta}
     \end{array}
    \right).
   \end{equation}
Henceforth, the equation \eqref{e1.1} changes to
  \begin{equation}\label{e7.5}
   (r^2\eta''+\varphi'')w_{11}-4r\eta'w_{1N}+2\eta w_{NN}+\Big(r^2\eta''+\varphi''-2r^2\eta'^2/\eta\Big)\Sigma_{i=2}^nw_{ii}=0,
  \end{equation}
where
  \begin{eqnarray*}
    w_{11}&=&-2\theta(2\eta)^{-n\theta}\psi^{-\theta-1}\Bigg(\eta''-\frac{2\eta'^2}{\eta}\Bigg)+4\theta(\theta+1)(2\eta)^{-n\theta}\psi^{-\theta-2}r^2\Bigg(\eta''-\frac{2\eta'^2}{\eta}\Bigg)^2\\
    w_{1N}&=&2n\theta^2(2\eta)^{-n\theta}\psi^{-\theta-1}\Bigg(\eta''-\frac{2\eta'^2}{\eta}\Bigg)\Bigg(r\frac{\eta'}{\eta}\Bigg)-2\theta(2\eta)^{-n\theta}\psi^{-\theta-1}\Bigg[r\Bigg(\eta''-\frac{2\eta'^2}{\eta}\Bigg)'\Bigg]\\
         &&+2\theta(\theta+1)(2\eta)^{-n\theta}\psi^{-\theta-2}\psi'\Bigg[r\Bigg(\eta''-\frac{2\eta'^2}{\eta}\Bigg)\Bigg]\\
    w_{ii}&=&-2\theta(2\eta)^{-n\theta}\psi^{-\theta-1}\Bigg(\eta''-\frac{2\eta'^2}{\eta}\Bigg), \ \ \ \ i=2,3,\cdots,n
  \end{eqnarray*}
and
   \begin{eqnarray*}
    w_{NN}&=&-n\theta(2\eta)^{-n\theta}\psi^{-\theta}\Bigg[\Bigg(\frac{\eta'}{\eta}\Bigg)'-n\theta\Bigg(\frac{\eta'}{\eta}\Bigg)^2\Bigg]+2n\theta^2(2\eta)^{-n\theta}\psi^{-\theta-1}\psi'\Bigg(\frac{\eta'}{\eta}\Bigg)\\
    &&+\theta(\theta+1)(2\eta)^{-n\theta}\psi^{-\theta-2}\psi'^2-\theta(2\eta)^{-n\theta}\psi^{-\theta-1}\psi''
   \end{eqnarray*}
hold for
   $$
    \psi\equiv r^2\Bigg(\eta''-\frac{2\eta'^2}{\eta}\Bigg)+\varphi'',\ \ \ \psi'\equiv\frac{\partial\psi}{\partial t}, \ \ \ \psi''\equiv\frac{\partial^2\psi}{\partial t^2}.
   $$
After substituting into \eqref{e7.5}, it yields that
   \begin{equation}\label{e7.6}
     Ar^4+Br^2+C=0,
   \end{equation}
in which $A, B, C$ are given by
   \begin{eqnarray*}
     A&\equiv&\Bigg(\eta''-\frac{2\eta'^2}{\eta}\Bigg)^2\Bigg\{4(\theta-n+1)\eta''+\Big[(2n^2-8n)\theta+(6n-4)\Big]\frac{\eta'^2}{\eta}\Bigg\}\\
      &&+(4n\theta-8\theta)\eta'\Bigg(\eta''-\frac{2\eta'^2}{\eta}\Bigg)\Bigg(\eta''-\frac{2\eta'^2}{\eta}\Bigg)'+2(\theta+1)\eta\Bigg[\Bigg(\eta''-\frac{2\eta'^2}{\eta}\Bigg)'\Bigg]^2\\
      &&-2\eta\Bigg(\eta''-\frac{2\eta'^2}{\eta}\Bigg)\Bigg(\eta''-\frac{2\eta'^2}{\eta}\Bigg)'',
   \end{eqnarray*}
together with
    \begin{eqnarray*}
      B&\equiv&\varphi''\Bigg\{4\Bigg[(\theta-2n+1)\eta''+(n^2\theta-2n\theta-2\theta+3n-3)\frac{\eta'^2}{\eta}\Bigg]\Bigg(\eta''-\frac{2\eta'^2}{\eta}\Bigg)\\
       &&+4(n\theta+2)\eta'\Bigg(\eta''-\frac{2\eta'^2}{\eta}\Bigg)'-2\eta\Bigg(\eta''-\frac{2\eta'^2}{\eta}\Bigg)''\Bigg\}\\
       &&+\varphi'''\Bigg\{4(n\theta-2\theta-2)\eta'\Bigg(\eta''-\frac{2\eta'^2}{\eta}\Bigg)+4(\theta+1)\eta\Bigg(\eta''-\frac{2\eta'^2}{\eta}\Bigg)'\Bigg\}\\
       &&+\varphi^{(4)}\Bigg\{-2\eta\Bigg(\eta''-\frac{2\eta'^2}{\eta}\Bigg)\Bigg\}
    \end{eqnarray*}
and
   \begin{eqnarray*}
      C&\equiv&\varphi''^2\Bigg\{-4n\eta''+2n(n\theta+3)\frac{\eta'^2}{\eta}\Bigg\}+4n\theta\eta'\varphi''\varphi'''\\
      &&+2(\theta+1)\eta\varphi'''^2-2\eta\varphi''\varphi^{(4)}.
    \end{eqnarray*}

\vspace{40pt}

\section{Euclidean complete affine maximal type hypersurfaces when $\theta\in(0,1/2)$}

When considering Euclidean complete solutions on half space ${\mathbb{R}}^n\times{\mathbb{R}}^+$, one can first choose $\eta=1$, a function fulfilling the positivity condition \eqref{e7.3}. Moreover, we have
  \begin{equation}\label{e8.1}
   \begin{cases}
    A\equiv0,\ \ \ \ B\equiv0,\\
    C=0\Leftrightarrow (\theta+1)\varphi'''^2-\varphi''\varphi^{(4)}=0.
   \end{cases}
  \end{equation}
Solving the fourth order equation of $\varphi$ yields that
   \begin{equation}\label{e8.2}
     \varphi(t)=C_1(t+C_2)^{2-\frac{1}{\theta}}+C_3t+C_4, \ \mbox{ or } \ \ C_5t^2+C_6t+C_7
   \end{equation}
for arbitrary constants $C_1-C_7$. In order that the hypersurface is convex, Euclidean complete and non-quadratic, we need only to choose $C_1>0$, $C_2=0$ in the former expression and let
    \begin{equation}\label{e8.3}
     2-\frac{1}{\theta}<0\Leftrightarrow0<\theta<1/2.
   \end{equation}
So, we achieve the following result.

\begin{theo}\label{t8.1}
  For each $N\geq2$ and $\theta\in(0,1/2)$, there exist Euclidean complete hypersurfaces on half space ${\mathbb{R}}^n\times{\mathbb{R}}^+$ which are not elliptic paraboloid. More precisely, they are given by
    \begin{equation}\label{e8.4}
      u(y,t)=C_1|y|^2+C_2t^{2-\frac{1}{\theta}}+C_3t+C_4, \ \ \forall t>0,
    \end{equation}
  where $C_1, C_2>0$ and $C_3, C_4\geq0$ are arbitrary constants.
\end{theo}

Secondly, if one chooses $\eta=t^{-1}, t>0$, the positivity condition \eqref{e7.3} for $\eta$ is also clearly true. Furthermore, we have
  \begin{equation}\label{e8.5}
   \begin{cases}
     A\equiv0, \ \ \ \ B\equiv0,\\
     C=0\Leftrightarrow n(n\theta-1)\varphi''^2-2n\theta t\varphi''\varphi'''\\
     \ \ \ \ \ \ \ \ \ \ \ \ \ \ \ +(\theta+1)t^2\varphi'''^2-t^2\varphi''\varphi^{(4)}=0.
   \end{cases}
  \end{equation}
Setting $\zeta\equiv t\frac{\varphi'''}{\varphi''}$, there holds
  \begin{equation}\label{e8.6}
    t\zeta'=\theta(\zeta+\beta_1)^2-\beta_2
  \end{equation}
for
   $$
     \beta_1\equiv\frac{1-2n\theta}{2\theta},\ \ \ \
     \beta_2\equiv\frac{1}{4\theta}.
   $$
Therefore, we get that
   \begin{equation}\label{e8.7}
     \zeta=\frac{n+(\theta^{-1}-n)\beta_3t}{1-\beta_3t}
   \end{equation}
for $\beta_3\in{\mathbb{R}}^-\cup\{0,-\infty\}$. After integrating $\zeta$ in variable $t$, it yields that
   \begin{equation}\label{e8.8}
    \varphi=\begin{cases}
      \beta_4t^{n+2}+\beta_5t+\beta_6, & \beta_3=0\\
      \beta_4\int^t_{(\beta_3/2)^{-1}}\int^s_{(\beta_3/2)^{-1}}\frac{r^n}{(r-\beta_3^{-1})^{\frac{1}{\theta}}}drds+\beta_5t+\beta_6, & \beta_3\in(-\infty,0)\\
      \beta_4t^{n+2-\frac{1}{\theta}}+\beta_5t+\beta_6, & \beta_3=-\infty
    \end{cases}
   \end{equation}
for constants $\beta_4>0$ and $\beta_5, \beta_6$. In order to be a convex large solution for $t>0$, we need to impose the condition $\beta_3=-\infty$ and
    \begin{equation}\label{e8.9}
       0<\theta<\frac{1}{n+2}
    \end{equation}
for $\varphi$. Summing up, we obtain a second result.

\begin{theo}\label{t8.2}
  For each $N\geq2$ and $\theta\in(0,1/(N+1))$, the hypersurfaces determined by
     \begin{equation}\label{e8.10}
       u(y,t)=C_1|y|^2t^{-1}+C_2t^{n+2-\frac{1}{\theta}}+C_3t+C_4, \ \ y\in{\mathbb{R}}^n, t>0
     \end{equation}
  are all Euclidean complete affine maximal type hypersurfaces on half space ${\mathbb{R}}^n\times{\mathbb{R}}^+$ for constants $C_1, C_2>0$ and $C_3, C_4\in{\mathbb{R}}$.
\end{theo}

\vspace{40pt}

\section{Trudinger-Wang type solutions for $\theta=1/2$ or $\theta=(N-1)/N$}

From now on, we try to look for the solutions of the form
   $$
    u(y,t)=\varphi(|y|)\eta(t).
   $$
The functions in this form have firstly been used by Trudinger-Wang \cite{TW1} to construct counter examples of Bernstein property of affine maximal hypersurfaces on dimension $N=10$, where they showed that the function $u$ given by
   $$
    n=9,\ \ \varphi(r)=r^9, \ \ \eta(t)=t^{-1}
   $$
satisfies \eqref{e1.1} on the set
   $$
     \Big\{(y,t)\in{\mathbb{R}}^n\times{\mathbb{R}}:\ y\not=0, t>0\Big\}
   $$
for $\theta=11/12$. Noting that the equation is invariant under rotation on ${\mathbb{R}}^N$, one may assume that $y=(r,0,\cdots,0)\in{\mathbb{R}}^n$ for simplicity. So, there holds
  \begin{eqnarray*}
    D^2u&=&\left(
     \begin{array}{ccccc}
       \varphi''\eta & 0 & \cdots & 0 & \varphi'\eta' \\
       0  &  \frac{\varphi'}{r}\eta & \cdots & 0 & 0\\
       \cdots & \cdots& \cdots & \cdots & \cdots\\
       0 & 0 & \cdots & \frac{\varphi'}{r}\eta & 0\\
      \varphi'\eta' & 0 & \cdots & 0 & \varphi\eta''
     \end{array}
    \right).
  \end{eqnarray*}
Hence, the inverse matrix of Hessian of $u$ is given by
  \begin{eqnarray*}
    [u^{ij}]&=&\left(
     \begin{array}{ccccc}
       \frac{\varphi\eta''}{\varphi\varphi''\eta\eta''-(\varphi'\eta')^2} & 0 & \cdots & 0 & -\frac{\varphi'\eta'}{\varphi\varphi''\eta\eta''-(\varphi'\eta')^2}\\
       0  &  \frac{r}{\varphi'\eta} & \cdots & 0 & 0\\
       \cdots & \cdots& \cdots & \cdots & \cdots\\
       0 & 0 & \cdots & \frac{r}{\varphi'\eta} & 0\\
      -\frac{\varphi'\eta'}{\varphi\varphi''\eta\eta''-(\varphi'\eta')^2} & 0 & \cdots & 0 & \frac{\varphi''\eta}{\varphi\varphi''\eta\eta''-(\varphi'\eta')^2}
     \end{array}
    \right)
  \end{eqnarray*}
and the determinant equals to
   \begin{equation}\label{e9.1}
     \det(D^2u)=\Big[\varphi\varphi''\eta\eta''-(\varphi'\eta')^2\Big]\Bigg(\frac{\varphi'\eta}{r}\Bigg)^{n-1}.
   \end{equation}
As a result,
   \begin{equation}\label{e9.2}
     w=\Big[\varphi\varphi''\eta\eta''-(\varphi'\eta')^2\Big]^{-\theta}\Bigg(\frac{\varphi'\eta}{r}\Bigg)^{-(n-1)\theta}
   \end{equation}
and the equation \eqref{e1.1} changes to
  \begin{equation}\label{e9.3}
   \varphi\eta''w_{11}-2\varphi'\eta'w_{1N}+\varphi''\eta w_{NN}+\frac{r}{\varphi'\eta}\Big[\varphi\varphi''\eta\eta''-(\varphi'\eta')^2\Big]\Sigma_{i=2}^nw_{ii}=0.
  \end{equation}
Taking $\varphi=e^{r^\alpha}, \eta=e^{t}$, then
  \begin{eqnarray}\nonumber\label{e9.4}
   w&=&(\varphi\varphi''-\varphi'^2)^{-\theta}\Bigg(\frac{\varphi'}{r}\Bigg)^{-(n-1)\theta}e^{-(n+1)\theta t}\\
   &=&\alpha^{-n\theta}(\alpha-1)^{-\theta}r^{\beta}e^{-(n+1)\theta(r^\alpha+t)}
  \end{eqnarray}
holds for $\beta\equiv-n\theta(\alpha-2)$. Henceforth,
  \begin{eqnarray*}
    \kappa w_{11}&=&\beta(\beta-1)r^{\beta-2}-(n+1)\theta\alpha(\alpha+2\beta-1)r^{\alpha+\beta-2}+(n+1)^2\theta^2\alpha^2r^{2\alpha+\beta-2}\\
    \kappa w_{1N}&=&-(n+1)\theta\beta r^{\beta-1}+(n+1)^2\theta^2\alpha r^{\alpha+\beta-1}\\
    \kappa w_{ii}&=&\beta r^{\beta-2}-(n+1)\theta\alpha r^{\alpha+\beta-2}\\
    \kappa w_{NN}&=&(n+1)^2\theta^2r^\beta
  \end{eqnarray*}
holds for $\kappa\equiv\alpha^{n\theta}(\alpha-1)^\theta e^{(n+1)\theta(r^\alpha+t)}$, and \eqref{e9.3} changes to
  \begin{equation}\label{e9.5}
    B_1r^{\beta-2}+B_2r^{\alpha+\beta-2}=0,
  \end{equation}
where
  \begin{eqnarray*}
    B_1&\equiv&\beta(\beta-1)+(n-1)(\alpha-1)\beta\\
    B_2&\equiv&-(n+1)\theta\alpha(\alpha+2\beta-1)+2(n+1)\theta\alpha\beta\\
     &&+\alpha(\alpha-1)(n+1)^2\theta^2-(n-1)(n+1)\alpha(\alpha-1)\theta.
  \end{eqnarray*}
So, we get a condition
  \begin{equation}\label{e9.6}
    \begin{cases}
       \beta\Big[(\beta-1)+(n-1)(\alpha-1)\Big]=0, \ \ \ \beta\equiv-n\theta(\alpha-2)\\
       (n+1)\theta\alpha(\alpha-1)\Big[(n+1)\theta-n\Big]=0,
    \end{cases}
  \end{equation}
and obtains the following result.

\begin{theo}\label{t9.1}
 Supposing that $N\geq2$ and $\theta=\frac{N-1}{N}$, there exists at least one non-quadratic entire affine maximal type hypersurface which is given by
   \begin{equation}\label{e9.7}
     u(y,t)=e^{|y|^\alpha+t}, \ \ y\in{\mathbb{R}}^n, t\in{\mathbb{R}}
   \end{equation}
for
   \begin{equation}\label{e9.8}
   \begin{cases}
    \alpha=2, &   N\geq2, \mbox{ or }\\
    \alpha=(N-1)(N-2), & N\geq3,
    \end{cases}
  \end{equation}
\end{theo}

It is remarkable that when $\alpha$ is taken to be $(N-1)(N-2), N\geq4$, the convexity of the solution was degenerated along the line $\{y=0, t\in{\mathbb{R}}\}$. As a corollary of Theorem \ref{t9.1}, we have the following result.

\begin{coro}\label{c9.1}
  Supposing that $N\geq2$ and $\theta=1/2$, there exist non-quadratic entire affine maximal type hypersurfaces which are given by
     \begin{equation}\label{e9.9}
       u(y,t)=e^{y_1^2+t}+\Sigma_{i,j=2}^{n} A^{ij}y_iy_j+\Sigma_{i=2}^nB^iy_i+C, \ \ \forall y\in{\mathbb{R}}^n, t\in{\mathbb{R}}
     \end{equation}
  for positive definite matrix $A\in{\mathbb{R}}^{(n-1)^2}$, vector $B\in{\mathbb{R}}^{n-1}$ and constant $C$.
\end{coro}

\vspace{40pt}

\section{Euclidean complete hypersurfaces in case of $\theta\in(1/2,(N-1)/N)$}

In the last section, we look for the solutions of the form
   $$
    u(y,t)=\Pi_{i=1}^ny_i^{-\alpha_i} e^t, \ \ y\in{\mathbb{R}}^n, y_i>0, t\in{\mathbb{R}}, \ \forall i=1,2,\cdots,n
   $$
for positive constants $\alpha_i, i=1,2,\cdots, n$.  A similar calculation shows that
  $$
    D^2u=\Pi_{i=1}^ny_i^{-\alpha_i} e^t\left(
     \begin{array}{ccccc}
       \alpha_1(\alpha_1+1)y_1^{-2} & \alpha_1\alpha_2y_1^{-1}y_2^{-1} & \cdots & \alpha_1\alpha_ny_1^{-1}y_n^{-1} & -\alpha_1y_1^{-1} \\
       \alpha_2\alpha_1y_1^{-1}y_2^{-1}  & \alpha_2(\alpha_2+1)y_2^{-2} & \cdots & \alpha_2\alpha_ny_2^{-1}y_n^{-1} &  -\alpha_2y_2^{-1}\\
       \cdots & \cdots & \cdots & \cdots & \cdots \\
      \alpha_n\alpha_1y_1^{-1}y_n^{-1} & \alpha_n\alpha_2y_2^{-1}y_n^{-1} & \cdots & \alpha_n(\alpha_n+1)y_n^{-2} & -\alpha_ny_n^{-1}\\
      -\alpha_1y_1^{-1} & -\alpha_2y_2^{-1} & \cdots & -\alpha_ny_n^{-1} & 1
     \end{array}
    \right).
  $$
So, the inverse matrix of Hessian matrix is given by
 $$
    [u^{ij}]=e^{-t}\Pi_{i=1}^n\alpha_i^{-1}y_i^{\alpha_i+2}\left(
     \begin{array}{ccccc}
       \Pi_{i\not=1}\alpha_iy_i^{-2} & 0 & \cdots & 0 & y_1\Pi_{i=1}^n\alpha_iy_i^{-2} \\
       0  & \Pi_{i\not=2}\alpha_iy_i^{-2} & \cdots & 0 &  y_2\Pi_{i=1}^n\alpha_iy_i^{-2}\\
       \cdots & \cdots & \cdots & \cdots & \cdots \\
      0 & 0 & \cdots & \Pi_{i\not=n}\alpha_iy_i^{-2} & y_n\Pi_{i=1}^n\alpha_iy_i^{-2}\\
      y_1\Pi_{i=1}^n\alpha_iy_i^{-2} & y_2\Pi_{i=1}^n\alpha_iy_i^{-2} & \cdots & y_n\Pi_{i=1}^n\alpha_iy_i^{-2} & (\Pi_{i=1}^n\alpha_iy_i^{-2})(\Sigma_{i=1}^n\alpha_i+1)
     \end{array}
    \right)
 $$
and the determinant equals to
   \begin{equation}\label{e10.1}
     \det(D^2u)=e^{Nt}\Pi_{i=1}^n\alpha_iy_i^{-N\alpha_i-2}.
   \end{equation}
Therefore,
   \begin{equation}\label{e10.2}
     w=e^{-N\theta t}\Pi_{i=1}^n\alpha_i^{-\theta}y_i^{(N\alpha_i+2)\theta}
   \end{equation}
and the equation \eqref{e1.1} changes to
  \begin{equation}\label{e10.3}
   \Sigma_{j=1}^n\Big(\Pi_{i\not=j}\alpha_iy_i^{-2}\Big)w_{jj}+2\Sigma_{j=1}^n\Big(y_j\Pi_{i=1}^n\alpha_iy_i^{-2}\Big)w_{jN}+\Big(\Pi_{i=1}^n\alpha_iy_i^{-2}\Big)\Big(\Sigma_{j=1}^n\alpha_j+1\Big)w_{NN}=0,
  \end{equation}
where
  \begin{eqnarray*}
    e^{N\theta t} w_{jj}&=&\Big[(N\alpha_j+2)\theta\Big]\Big[(N\alpha_j+2)\theta-1\Big]y_j^{-2}\Pi_{i=1}^n\alpha_i^{-\theta}y_i^{(N\alpha_i+2)\theta}\\
    e^{N\theta t} w_{jN}&=&-N\theta\Big[(N\alpha_j+2)\theta\Big]y_j^{-1}\Pi_{i=1}^n\alpha_i^{-\theta}y_i^{(N\alpha_i+2)\theta}\\
    e^{N\theta t} w_{NN}&=& N^2\theta^2\Pi_{i=1}^n\alpha_i^{-\theta}y_i^{(N\alpha_i+2)\theta}.
  \end{eqnarray*}
As a result, \eqref{e10.3} changes to
  \begin{equation}\label{e10.4}
   \Sigma_{j=1}^n\alpha_j^{-1}\Big[(N\alpha_j+2)\theta\Big]\Big[(N\alpha_j+2)\theta-1\Big]-2N\theta^2\Sigma_{j=1}^n(N\alpha_j+2)+N^2\theta^2\Big(\Sigma_{j=1}^n\alpha_j+1\Big)=0,
  \end{equation}
or equivalently
   \begin{equation}\label{e10.5}
    \theta=\frac{2\Sigma_{j=1}^n\alpha_j^{-1}+N(N-1)}{4\Sigma_{j=1}^n\alpha_j^{-1}+N^2}.
   \end{equation}
So, we obtain the following result.

\begin{theo}\label{t10.1}
 Supposing that $N\geq3$ and
   \begin{equation}\label{e10.6}
     \theta\in(1/2, (N-1)/N),
   \end{equation}
there exist Euclidean complete affine maximal type hypersurfaces on
  $$
   {\mathbb{H}}\equiv\Big\{(y,t)\in{\mathbb{R}}^n\times{\mathbb{R}}\big|\ y_i>0, \ \forall i=1,2,\cdots,n\Big\},
  $$
which are given by
   \begin{equation}\label{e10.7}
     u(y,t)=\Pi_{i=1}^ny_i^{-\alpha_i}e^t, \ \ y\in{\mathbb{R}}^n, y_i>0, t\in{\mathbb{R}}, \ \ \forall i=1,2,\cdots, n
   \end{equation}
for some positive constants $\alpha_i, i=1,2,\cdots, n$.
\end{theo}

\noindent\textbf{Proof.} For $N\geq3$ and $\theta\in(1/2,(N-1)/N)$, one can choose positive constants $\alpha_j, j=1,2,\cdots, n$ such that
     $$
      \frac{2\Sigma_{j=1}^n\alpha_j^{-1}+N(N-1)}{4\Sigma_{j=1}^n\alpha_j^{-1}+N^2}=\theta
     $$
by continuity. So, the above construction gives the desired Euclidean complete affine maximal type hypersurfaces on ${\mathbb{H}}$. The proof was done. $\Box$\\

When $\theta$ belongs to the range $(1/2,(N-1)/N)$ for dimension $N\geq3$, one may also try the solutions of the form
   $$
    u(y)=\Pi_{i=1}^Ny_i^{-\alpha_i}, \ \ y\in{\mathbb{R}}^n, y_i>0, \ \forall i=1,2,\cdots,N
   $$
for positive constants $\alpha_i, i=1,2,\cdots, n$.  A similar calculation shows that
  $$
    D^2u=\Pi_{i=1}^Ny_i^{-\alpha_i}\left(
     \begin{array}{cccc}
       \alpha_1(\alpha_1+1)y_1^{-2} & \alpha_1\alpha_2y_1^{-1}y_2^{-1} & \cdots & \alpha_1\alpha_Ny_1^{-1}y_N^{-1}\\
       \alpha_2\alpha_1y_1^{-1}y_2^{-1}  & \alpha_2(\alpha_2+1)y_2^{-2} & \cdots & \alpha_2\alpha_Ny_2^{-1}y_N^{-1}\\
       \cdots & \cdots & \cdots & \cdots \\
      \alpha_N\alpha_1y_1^{-1}y_N^{-1} & \alpha_N\alpha_2y_2^{-1}y_N^{-1} & \cdots & \alpha_N(\alpha_N+1)y_N^{-2}
     \end{array}
    \right).
  $$
So, the inverse matrix is given by
 $$
    [u^{ij}]=\beta^{-1}\Pi_{i=1}^N\alpha_i^{-1}y_i^{\alpha_i+2}\left(
     \begin{array}{cccc}
       \beta_1\Pi_{i\not=1}\alpha_iy_i^{-2} & -y_1y_2\Pi_{i=1}^N\alpha_iy_i^{-2} & \cdots & -y_1y_N\Pi_{i=1}^N\alpha_iy_i^{-2}\\
       -y_1y_2\Pi_{i=1}^N\alpha_iy_i^{-2}  & \beta_2\Pi_{i\not=2}\alpha_iy_i^{-2} & \cdots & -y_2y_N\Pi_{i=1}^N\alpha_iy_i^{-2}\\
       \cdots & \cdots & \cdots & \cdots \\
      -y_1y_N\Pi_{i=1}^N\alpha_iy_i^{-2} &  -y_2y_N\Pi_{i=1}^N\alpha_iy_i^{-2} & \cdots & \beta_N\Pi_{i\not=N}\alpha_iy_i^{-2}
     \end{array}
    \right)
 $$
and determinant equals to
   \begin{equation}\label{e10.8}
     \det(D^2u)=\beta\Pi_{i=1}^N\alpha_iy_i^{-N\alpha_i-2},
   \end{equation}
where
   $$
   \beta\equiv\Sigma_{i=1}^N\alpha_i+1, \ \ \  \beta_j\equiv \Sigma_{i\not=j}\alpha_i+1
   $$
Therefore,
   \begin{equation}\label{e10.9}
     w=\beta^{-\theta}\Pi_{i=1}^N\alpha_i^{-\theta}y_i^{(N\alpha_i+2)\theta}
   \end{equation}
and the equation \eqref{e1.1} changes to
  \begin{equation}\label{e10.10}
   \Sigma_{j=1}^N\alpha_j^{-1}\beta_jy_j^2w_{jj}-2\Sigma_{j\not=k}y_jy_kw_{jk}=0,
  \end{equation}
for
  \begin{eqnarray*}
    w_{jj}&=&\beta^{-\theta}\Big[(N\alpha_j+2)\theta\Big]\Big[(N\alpha_j+2)\theta-1\Big]y_j^{-2}\Pi_{i=1}^N\alpha_i^{-\theta}y_i^{(N\alpha_i+2)\theta}\\
    w_{jk}&=&\beta^{-\theta}\Big[(N\alpha_j+2)\theta\Big]\Big[(N\alpha_k+2)\theta\Big]y_j^{-1}y_k^{-1}\Pi_{i=1}^N\alpha_i^{-\theta}y_i^{(N\alpha_i+2)\theta}
  \end{eqnarray*}
As a result, \eqref{e10.10} can be simplified by
  \begin{equation}\label{e10.11}
   \Sigma_{j=1}^N\alpha_j^{-1}\beta_j\Big[(N\alpha_j+2)\theta\Big]\Big[(N\alpha_j+2)\theta-1\Big]-2\Sigma_{j<k}\Big[(N\alpha_j+2)\theta\Big]\Big[(N\alpha_k+2)\theta\Big]=0
  \end{equation}
or equivalently
   \begin{eqnarray}\nonumber\label{e10.12}
    \theta&=&\frac{\Sigma_{j=1}^N\alpha_j^{-1}\beta_j(N\alpha_j+2)}{\Sigma_{j=1}^N\alpha_j^{-1}\beta_j(N\alpha_j+2)^2-2\Sigma_{j<k}(N\alpha_j+2)(N\alpha_k+2)}\\
    &=&\frac{\beta\Sigma_{j=1}^N\alpha_j^{-1}(N\alpha_j+2)-\Sigma_{j=1}^N(N\alpha_j+2)}{\beta\Sigma_{j=1}^N\alpha_j^{-1}(N\alpha_j+2)^2-\Big[N\Sigma_{j=1}^N\alpha_j+2N\Big]^2}\\ \nonumber
    &=&\frac{2\beta\Sigma_{j=1}^N\alpha_j^{-1}+(N^2-N)\beta-N}{4\beta\Sigma_{j=1}^N\alpha_j^{-1}+N^2\beta-N^2}.
   \end{eqnarray}
Setting $\gamma=\Sigma_{j=1}^N\alpha_j^{-1}+1$, the last formula is equivalent to
  \begin{equation}\label{e10.13}
   (4\theta-2)\beta\gamma+(N^2-4)\Bigg(\theta-\frac{N+1}{N+2}\Bigg)\beta=N(N\theta-1).
  \end{equation}
Now, for $\theta\in(1/2,(N+1)/(N+2))$, we want to find out the range of the value of the function
  $$
   F(\alpha)\equiv(4\theta-2)\beta\gamma+(N^2-4)\Bigg(\theta-\frac{N+1}{N+2}\Bigg)\beta
  $$
on the domain
  $$
   {\mathbb{Q}}\equiv\Bigg\{\alpha=(\alpha_1,\cdots,\alpha_N)\in{\mathbb{R}}^N\Big|\ \alpha_i>0, \ \forall i=1,2,\cdots,N\Bigg\}
  $$
Supposing first that $\theta>\frac{N-1}{N}$, then
  $$
   \begin{cases}
     F_{\alpha_1}=(4\theta-2)\gamma-(4\theta-2)\beta\alpha_1^{-2}+(N^2-4)\Big(\theta-\frac{N+1}{N+2}\Big)=0\\
     F_{\alpha_2}=(4\theta-2)\gamma-(4\theta-2)\beta\alpha_2^{-2}+(N^2-4)\Big(\theta-\frac{N+1}{N+2}\Big)=0\\
    \  \ \ \ \ \ \ \ \ \ \ \  \cdots \cdots \cdots \cdots \cdots \cdots \cdots \cdots \cdots \cdots \ \ \ \ \  \\
     F_{\alpha_N}=(4\theta-2)\gamma-(4\theta-2)\beta\alpha_N^{-2}+(N^2-4)\Big(\theta-\frac{N+1}{N+2}\Big)=0
   \end{cases}
  $$
gives the unique critical point of $F$ on ${\mathbb{Q}}$ by
   $$
    \alpha_1=\alpha_2=\cdots=\alpha_N\equiv x,
   $$
where $x$ is given by
  $$
   x=\frac{1}{N}\sqrt{\frac{4\theta-2}{\theta-\frac{N-1}{N}}}.
  $$
Noting that for $\theta>(N-1)/N$, there holds
   $$
    (4\theta-2)+(N^2-4)\Bigg(\theta-\frac{N+1}{N+2}\Bigg)>0.
   $$
So, the function $F$ tends to infinity as $||\alpha||\to\infty, \alpha\in{\mathbb{Q}}$. This means that
  \begin{equation}\label{e10.14}
   \inf_{\alpha\in{\mathbb{Q}}} F(\alpha)=F(x)=5N^2\theta-(3N^2-N)+2N^2\sqrt{(4\theta-2)\Bigg(\theta-\frac{N-1}{N}\Bigg)}.
  \end{equation}
As a result, in order that there exists at least one $\alpha\in{\mathbb{Q}}$ such that \eqref{e10.13} is fulfilled for given $\theta\in((N-1)/N,(N+1)/(N+2))$, we need
  \begin{eqnarray*}
   && \inf_{\alpha\in{\mathbb{Q}}} F(\alpha)<N(N\theta-1)\\
   &\Leftrightarrow&4N^2\theta-(3N^2-2N)+2N^2\sqrt{(4\theta-2)(\theta-\frac{N-1}{N})}<0\\
   &\Rightarrow&\theta<\frac{3}{4}-\frac{1}{2N}.
  \end{eqnarray*}
which contradicts with our assumption $\theta>\frac{N-1}{N}$ for $N\geq3$. Now, we turn to consider the case $\theta\in(1/2,(N-1)/N)$. Noting that
   $$
    (4\theta-2)+(N^2-4)\Bigg(\theta-\frac{N+1}{N+2}\Bigg)=N^2\theta-N(N-1)<0
   $$
in this case, if one lets $\alpha_1, \alpha_2, \cdots, \alpha_N$ become all large, there hold
   $$
    \gamma\to 1, \ \ \beta\to+\infty
   $$
and hence $\inf_{{\mathbb{Q}}}F(\alpha)=-\infty$. On another hand, if one lets $\alpha_1, \alpha_2, \cdots, \alpha_N$ become all small, there hold
   $$
    \gamma\to +\infty, \ \ \beta\to1
   $$
and thus $\sup_{{\mathbb{Q}}}F(\alpha)=+\infty$. So, we reach the following result.

\begin{theo}\label{t10.2}
 Supposing that $N\geq3$ and
   \begin{equation}\label{e10.15}
     \theta\in(1/2,(N-1)/N),
   \end{equation}
 there exist Euclidean complete affine maximal type hypersurfaces on
    $$
     \Big\{y\in{\mathbb{R}}^N\big|\ y_i>0, \ \forall i=1,2,\cdots, N\Big\}
    $$
 which are given by
   \begin{equation}\label{e10.16}
    u(y)=\Pi_{i=1}^ny_i^{-\alpha_i}, \ \ y\in{\mathbb{R}}^n, y_i>0, \ \forall i=1,2,\cdots,N
   \end{equation}
 for some positive constants $\alpha_i, i=1,2,\cdots,N$.
\end{theo}

Now, the conclusion of Theorem C follows from Theorem \ref{t8.1} for $\theta\in(0,1/2)$, Theorem \ref{t9.1} together with Corollary \ref{c9.1} for $\theta=1/2 \mbox{ or } (N-1)/N$, and Theorem \ref{t10.1} or Theorem \ref{t10.2} for $\theta\in(1/2,(N-1)/N)$. $\Box$\\

\vspace{40pt}

\section*{Acknowledgments}

The author would like to express his deepest gratitude to Professors Xi-Ping Zhu, Kai-Seng Chou, Xu-Jia Wang and Neil Trudinger for their constant encouragements and warm-hearted helps. This paper is also dedicated to the memory of Professor Dong-Gao Deng.\\

\vspace{10pt}



\end{document}